\magnification = \magstephalf \hsize=14truecm \hoffset=1truecm \parskip=5pt.
\input amssym.def

\def\Rb{{\overline{R}}}
\def\pb{{\overline{p}}}
\def\HH{{\cal H}}
\def\IR{{\rm I}\kern-.2em{\rm R}}
\baselineskip=18pt
\font\eightit=cmti8
\font\eightpoint=cmr8

\def\lanbox{\hbox{$\, \vrule height 0.25cm width 0.25cm depth
0.01cm \,$}}

\vskip1.4cm \noindent

\hbox to \hsize{\hfil\eightit 31aug00} \vglue1truein

\centerline{\bf A direct proof of a theorem of Blaschke and Lebesgue }
\medskip
{\baselineskip=2.5ex
\vfootnote{}{\eightpoint
\noindent\copyright 2000 by the author.
Reproduction of this article, in its entirety, by any means is permitted
for non--commercial purposes.}}

{\baselineskip = 12pt

\centerline{\bf Evans M. Harrell II} 
\footnote {$^{*}$}{harrell@math.gatech.edu, School of Mathematics,
Georgia Tech, Atlanta, GA 30332--0160, USA.} 

\smallskip

\vskip .7 true in

\centerline{\bf Abstract}

\vskip .7 true in

\smallskip
The Blaschke--Lebesgue Theorem
states that among all planar convex domains of given constant width
$B$ the Reuleaux triangle has minimal area.  It is the purpose of the
present note to give a direct proof of this theorem by
analyzing the underlying variational problem.  
The advantages of the
proof are that it shows uniqueness 
(modulo rigid deformations such as rotation
and translation) and leads analytically to
the shape of the area--minimizing
domain.  
Most previous proofs have relied on foreknowledge 
of the minimizing domain.
Key parts of the analysis extend to the
higher--dimensional situation, 
where the convex body of given constant width and
minimal volume is unknown. 

Mathematics subject classification 52A10, 52A15, 52A38, 49Q10

\noindent
{\vglue 0.5cm}
\medskip

\noindent{\bf I. Introduction.}
\medskip

A convex body in $\IR^d$ is said to have constant width $B$ 
if any two distinct parallel planes tangent to 
its boundary are separated by a distance $B.$
For $d=2$ such bodies are often called {\it orbiforms},
and for $d=3$ they are called {\it spheroforms}. A well-known 
example is the
{\it Reuleaux triangle}, whose boundary consists 
of three equally long circular arcs
with curvature $1/B$. The arcs meet at the corners of an equilateral
triangle.  Reuleaux polygons with any odd number of sides
likewise enjoy the property of constant width

It has long been known that among all two-dimensional
convex bodies of constant width, the Reuleaux triangle
has the smallest area.
W. Blaschke [Bla15] and H. Lebesgue [Leb14, Leb21]
were the first to show this,
and the succeeding
decades have seen several other works on the the problem of
minimizing the area or volume of an object given a constant
width; see [Fuj27-31, BoFe34, Egg52, Bes63, ChGR83, Web94, 
Gha96, and CCG96].   Objects of constant width have several 
practical uses, and have been entertainingly discussed in
[Fey89, Kaw98].  For instance, 
coins are sometimes designed with such shapes, because 
constant width allows their use in vending machines.
 
Ths disadvantage of the arguments of Blaschke and Lebesgue and 
most subsequent proofs of the Blaschke-Lebesgue theorem is that they are not suf ficiently analytic to derive the minimality of the 
Reuleaux triangle without prior knowledge.  No doubt this is a
one of the reasons that the higher-dimensional analogue of the problem
has remained open:  What body (or bodies) of constant width in 
three or more dimensions has the smallest volume?  

It is my purpose here to prove the Blaschke-Lebesgue theorem in
a directly analytic way, and to frame the problem 
in higher dimensions
as a step toward answering the question just posed.  

Two previous attempts to give 
analytical proofs can be cited.  
Fujiwara [Fuj27-31] expressed the area in terms of $r(\theta)$ and  showed throu gh a lengthy 
calculation that in general
the area of an orbiform exceeds that of the Reuleaux triangle.  His
proof gives
little indication how to find the optimal geometry from
first principles.  
More recently Ghandehari [Gha96] gave an argument 
via optimal control theory and Pontryagin's maximum principle,
which resembles the one of this article in a few respects.

\medskip
\noindent{\bf II.  On the Minimal Volume of a Convex Body 
of Constant Width.}
\medskip

A body $D$ of constant width is strictly convex, and 
therefore $\partial \Omega$ may be expressed as a 
continuous image 
of the sphere $S^{d-1}$ via the mapping 
$\Gamma({\bf \omega})$ which
associates to any unit vector ${\bf \omega}$ the point of 
$\partial D$ farthest in the 
direction ${\bf \omega}.$  
(At smooth points of $\partial D, \Gamma$ is
the inverse of the Gauss map.)

If ${\bf x}$ denotes a point on the boundary, then
the {\it support function} of $D$ will be defined
in the usual manner as
$p({\bf \omega}):={\bf x} \cdot {\bf \omega}.$  
Notice that $p({\bf \omega})$ is the distance from the origin of a plane in contact with $\partial D,$ provided that the origin is
within $D,$ which may always be assumed.  Once
$p({\bf \omega})$ is known, one can reconstruct the convex set as
the envelope of its supporting planes.  Choosing the independent variable as 
${\bf \omega}$ will be convenient for several reasons, 
among them the simple form of the formula for the 
width of $D$:

$$p({\bf \omega})+p({\bf \omega}^{\bf a})=B, \eqno(1)$$

\noindent
where ${\bf \omega}^{\bf a}$ designates the point on 
$S^1$ antipodal to ${\bf \omega}$: for $S^1$ one could 
identify $\omega$ as the usual angular variable and
write ${\omega}^{\bf a} = {\omega} \pm \pi,$
but dimension--independent notation will be preferred 
as far as possible.

A simple exercise using the divergence theorem shows that the
volume can be written in terms of the support function:

$$Vol ={1\over d}\int p({\bf \omega}) dS
= {1\over d}\int p({\bf \omega}) 
{d\omega \over {\sum_{j} \kappa_j}} 
$$

In this formula $\kappa_j$ are the principal curvatures of
$\partial D.$  Here and elsewhere, it will be more 
convenient to express 
quantities in terms of the radii of curvature 
$R_j := {1 \over \kappa_j}$
(or zero,
at non-smooth points of $\partial D$).  Hence

$$
Vol =\left|{1\over d}\left\langle p, 
\prod^{d-1}_{j=1} R_j\right\rangle_{S^{d-1}} \right|. \eqno(2)$$

\noindent
The brackets here designate the inner product on $L^2(S^{d-1}).$
The set-up described to this stage is classical; for instance see [Bla49].

The question under consideration is:

\noindent
{\bf Problem 1}: Determine the minimal volume of a convex body of fixed
width $B.$

This problem will be recast with the benefit of several observations, beginning  with a useful formula,
which results from a direct calculation:

$$\nabla^2_{S^{d-1}} p = \sum^{d-1}_{j=1} R_j-(d-1) p. \eqno(3)$$

\noindent
Equation (3) was known to Weingarten [Wei1884] in the
nineteenth century.  Together with $(1)$ it implies that

$$\sum_j R_j({\bf \omega})+\sum_j R_j({\bf \omega}^{\bf a})=(d-1) B. \eqno(4)$$

Observe that $d-1$ is the second eigenvalue of $-\nabla^2_{S^{d-1}},$ 
so the 
differential equation $(3)$ is not uniquely solvable.  
However, according to the 
Fredholm alternative theorem,
it is uniquely solvable with a reduced resolvent
$G: \HH_1 \hookleftarrow$, where 

$$\HH_1 := L^2(S^{d-1}) \ominus \hbox{ span}[Y^m_1],$$
and $Y^m_1$ are the spherical harmonics [M\"ul66] such that
$$-\nabla^2_{S^{d-1}} Y^m_1=(d-1)Y^m_1$$

\noindent
(If $d=2$, then $Y^m_1=\sin(\omega)$ and $\cos\omega$.  The notation 
$Y^m_\ell$ will be used in any dimension.) 

Now, $G$ is a bounded smoothing, operator.  That is,
$$Vol ={1\over d}\left\langle 
G\left[\sum\limits_{j} R_j\right], \prod\limits_{j} R_j\right\rangle_{S^{d-1}}
\eqno(5)$$
is a bounded quadratic form on $L^2(S^{d-1}),$ and the operator
G maps $L^2(S^{d-1})$ into 
$W^2(S^{d-1}) \cap \HH_1.$
Moreover, the condition of orthogonality to the span of the 
$ Y^m_1 $
is quite natural geometrically.  For the support function, this 
restriction means that the centroid has been moved to the
center of the coordinate system.  Any given set of 
nonzero coefficients 
of $ Y^m_1 $ could be specified, and this 
would merely correspond to
rigidly displacing the body $\Omega$ by a fixed vector with respect
to the centroid.  On the other hand, the condition that
$\prod^{d-1}_{j=1} R_j$ be orthogonal to $ Y^m_1 $ 
is necessary for $\partial \Omega$ to be a closed boundary:  If 
$d=2,$ it is the condition that the curve $\partial \Omega$ be closed.
If $d=3,$ this condition is necessary and essentially sufficient for 
the Gauss curvature to determine an immersed closed surface 
$\partial \Omega$ (uniquely up to rigid motions) [Min67, p. 130].

When $d=2$, there is only one curvature defined on
the boundary, and $V$ becomes a symmetric quadratic
form
$$Vol={1\over 2}\langle G[R], R\rangle_{S^1}\eqno(6)$$

When $d=3$, a theorem of Blaschke [see ChGr83, p 66] states that for
objects of constant width $B,$ the volume and surface--area $S$ are
related by 
$$Vol = {{B S} \over 2} - {{\pi B^3} \over 3}.$$
It follows that the minimizers
of the volume functional are identical to the minimizers of the 
surface--area functional, which may be written as a symmetric
quadratic form in $R := \sum\limits_{j} R_j$:

$$\Phi_1[R] := {1\over d}\langle G[R], R\rangle_{S^{d-1}}.\eqno(7)$$

\noindent
Recall that the support function enters through $G[R] = p.$  
The functional $(7)$ will be considered here as the objective 
in any dimension, although its interpretation is not immediate 
when $d > 3.$

With this notation, $(4)$ is written:

$$R({\bf \omega})+R({\bf \omega}^{\bf a})= (d-1)B, \eqno(8)$$

\noindent
This implies that admissible $R$ must satisfy 

$$0 \leq R({\bf \omega}) \leq  (d-1)B, \eqno(9)$$

\noindent
and the averages of $R$ and $p$ are both determined:  
It follows from $(8)$ and $(1)$ that

$$R_{ave} = {{(d-1) B} \over 2}, p_{ave} = B/2. \eqno(10)$$

Hence a simplification is achieved 
by subtracting the averages
of $R$ and $p$, so $\Rb := R - {{(d-1) B} \over 2}$ and 
$\pb := p - {B \over 2}.$
In these terms, just as
$p = G[R], \pb = G[\Rb]$.  
There results an alternative to Problem 1:  

\noindent
{\bf Problem 2}: Minimize the functional
$$\Phi[\Rb] := \langle G[\Rb], \Rb \rangle_{S^{d-1}} \eqno(11)$$

\noindent
for $\Rb \in \HH := \left\{ f \in L^2(S^{d-1}): f \perp 
\left\{Y_1^m \right\}, f({\bf \omega^a}) = - f({\bf \omega}),
| f({\bf \omega})| \leq {(d-1)B \over 2} \right\}.$

\noindent
{\bf Remarks}.  

\noindent
1.  Functions in $\HH$ are orthogonal to the lowest two 
eigenspaces of $-\nabla^2$.  It follows that $\Phi$
is a negative definite quadratic form on $\HH.$  In
particular, the function corresponding to the ball,
$\Rb = 0,$ maximizes
$\Phi.$  Because of the convexity of $\Phi,$
the minimizers are extremals of
$\HH.$  This statement is made somewhat
more precise in Theorem 1, below.

\noindent
2.  When $d=2,$ minimizing $\Phi$ on $\HH$ 
is equivalent to finding the convex region of smallest area
for a given $B.$  When 
$d=3,$ the theorem of Blaschke alluded to above ensures that
minimizing $\Phi$ is equivalent to minimizing the volume
functional, but some elements of $\HH$ may not correspond to 
embedded convex bodies.  Hence Problem 2 is fully equivalent to
Problem 1 only for $d=2.$

Now, the derivative of $\Phi$ with respect to
the variation $\Rb \to \Rb + \delta \zeta$ 
is simply 
$${d\Phi \over d\delta}\Big|_0 =2\,\left< G[\Rb],\zeta
\right>_{S^{d-1}}= 2\,\left< \pb, \zeta\right>.
 \eqno(12)$$

It is now possible to conclude:

\noindent
{\bf Theorem 1}. 
{\it Minimizers of Problem 2 exist, and
every minimizing $\Rb$
has the properties that}
\item{(a)} $\mu\left\{{\bf \omega}: \pb > 0,
|\Rb| < {(d-1)B \over 2} \right\}=0$
\item{(b)} $\pb > 0 \Rightarrow 
\Rb= -{(d-1)B \over 2}$, $\pb < 0
\Rightarrow \Rb= {(d-1)B \over 2}$
{\it a.e.}

\noindent
{\bf Proof}.  
The existence of a minimizer follows in a standard 
way from the compactness of the operator $G,$ considered as 
an operator on the Hilbert space
$$\left\{ f \in L^2(S^{d-1}): f \perp 
\left\{Y_1^m, 1 \right\} \right\}.$$  
(Minimizers are non--unique 
at least by rotation.)  

Consider now admissible variations for $\Phi$, 
normalizing $B$ temporarily for convenience so that
${(d-1)B \over 2}=1,$ and thus $-1 \leq \Rb \leq 1.$ 

Suppose that for some minimizing $\Rb$ and some
$\epsilon > 0,$
the set 
$ S_\epsilon := \left\{{\bf \omega}: \pb > 0,
|\Rb| \leq {(d-1)B \over 2}-\epsilon\right\}$
is of positive measure.  Then 
the antipodal set $S^{\bf a}_{\epsilon}$
is also of positive measure,
and any variation
$\zeta$ supported in $S_\epsilon$
must be extended to $S^{\bf a}_\epsilon$ antisymmetrically
by $\zeta ({\bf \omega}^{\bf a})=-\zeta ({\bf \omega})$.
Observe here that it is unnecessary to restrict $\zeta$ to be
orthogonal to $Y_1^m,$ as any such component 
is orthogonal to $\pb$ and hence will not
contribute to (12).

Now let $\zeta$ run through a basis for 
$L^2[S_{\epsilon}]\ominus\chi_  {S_\epsilon}$ consisting
of bounded functions $\zeta_n$.  (Boundedness,
together with $\epsilon > 0$ ensures admissibility.
The case $\zeta$ proportional to $\chi_{S_\epsilon}$
will be considered separately below.)
From $(11),$ 
with $\pb\left( {\bf \omega} \right) := p[\Rb]$ 
the first variation $(12)$ is proportional to

$$\eqalign{&\left\langle {\pb,\zeta } \right\rangle =\int\limits_{S_\varepsilon    } {\pb\left( {\bf \omega}  \right)\zeta \left( {\bf \omega}  \right)d\omega +}  \int\limits_{S_\varepsilon ^{\bf a}}
{\pb\left( {\bf \omega}  \right)\zeta \left( {\bf \omega}  \right)d\omega }\cr
  &\        =\int\limits_{S_\varepsilon } {\pb\left( {\bf \omega}  \right)\zeta    \left( {\bf \omega}  \right)d\omega} + \int\limits_{S_\varepsilon } {\left( {- \pb\left( {\bf \omega}  \right)}
\right) \left(-\zeta \left( {\bf \omega}  \right) \right) d\omega }\cr
  &\        
=2\int\limits_{S_\varepsilon } {\pb\left( {\bf \omega}  \right)\zeta \left( {\bf    \omega}  \right)d\omega .}\cr}$$

\noindent
Optimality implies that this vanishes and 
hence that $\pb =$ constant a.e.\ on $S_\epsilon$.  

Next consider $(11)$ subjected to the variation
$\zeta=-\chi_{S_\epsilon} + \chi_{S_\epsilon^{\bf a}}$:  

\noindent
If
$\mu(S_\epsilon)>0$, then
$${dV\over d\delta}=-\int_{S_\epsilon}\pb  + 
\int_{S_\epsilon^{\bf a}} \pb < 0,  \eqno(13)$$
which contradicts optimality.  This concludes the proof of (a).

For (b), observe from (a) that either
$\pb = 0$
a.e., which corresponds to the sphere, i.e., 
the maximizing shape,
or else there is a set of positive measure for which 
$p  > 0$ and $\Rb=-1$
or $+1$. But if $R=+1$, then the variation leading
to $(13)$ is still admissible for
$\delta \ge 0,$ so $(13)$ yields a contradiction.
Similarly for $p < 0$ if $R=-1.$
\hfil \lanbox

\noindent
{\bf Corollary 2}.  (The Blaschke--Lebesgue theorem.)  
{\it Among all two--dimensional convex 
regions of a given constant width $B$, the Reuleaux 
triangle has the smallest area.}

\noindent
{\bf Proof}.  Here $\omega$ is treated as the angular variable 
for $S^1$, and it will be assumed that $B=1.$
As the circle is not the minimizer, statement (b) of
Theorem 1, implies that $m := \max \pb  > 0$.  
By performing a rotation, it may be assumed that 
$\pb(0) = m,$ and
by continuity there is an interval around 
$0$ such that, when rewritten in terms of 
$\pb$ and specialized to one variable, (3) becomes
$$\pb''=-\pb - {1 \over 2}, \eqno(14)$$
yielding
$$\pb =\left( m + {1 \over 2}\right) \cos(\omega) - {1 \over 2} \eqno(15)$$
on that interval.  The end points of the interval correspond to 
$\pb =0$, i.e., 
$\omega= \pm \arccos {1\over {2m+1}}
=: \pm \alpha.$  At these points, $\pb' \ne 0$. 
Since standard regularity theory implies that
$\pb  \in AC^1$, $\pm \alpha$ cannot 
abut an interval on
which  $\pb =0.$  The only possibility is that
$\pb$ becomes negative and on the next interval
the differential equation (14) produces a solution
antisymmetric about $\alpha,$ i.e.,
$$\pb = - \left( m + {1 \over 2}\right) \cos(2 \alpha -\omega) + {1 \over 2}. 
\eqno(16)$$

The function $\pb$ switches between the two forms (15) and (16).


\noindent
The support function is also
subject to periodicity ($\omega + 2 \pi \cong \omega$) and
antisymmetry ($\pb(\omega+\pi) = -\pb(\omega)$).
The only candidates for optimality correspond to the
odd--sided Reuleaux polygons with
$R=1$.  An elementary calculation shows that the area of 
any such figure of given width is an increasing function of the
number of sides.
\hfill \lanbox

{\bf Concluding Remarks}.  

Two specific barriers have so far prevented the extension 
of this analysis to higher dimensions.  One of these is connected 
with the 
ability to extend solutions of ordinary differential equations uniquely
across a boundary; this needs to be replaced by a PDE analysis.

The second barrier is geometrical, elucidating
the nature of the set in the the $R_1$--$R_2$ plane 
which corresponds to convex bodies of constant width
$B$ in terms of $R:= R_1+R_2.$

\medskip
\noindent{\bf {Acknowledgments}}
\medskip

I am very grateful to W. Gangbo
and B. Kawohl, without whom 
this article would not have existed.  They
informed me of the problem and 
provided numerous useful comments and references.
They also passed on comments and references to me
from M. Belloni and E. Heil to whom I am thus indirectly 
indebted.  Finally, I enjoyed the hospitality of J. Fleckinger
at CEREMATH in Toulouse while some of this work was done.

\medskip
\noindent{\bf {References}}
\medskip

\item{[Bes63]} 
A.S. Besicovich, Minimum area of a set of constant width,
{\it Proc. Symp. Pure Math.} {\bf 7} (1963)13--14.

\item{[Bla15]} 
W. Blaschke, 
Konvexe Bereiche gegebener konstanter Breite
und kleinsten Inhalts,
{\it Math.  Ann.} {\bf 76} (1915)504--513.

\item{[Bla49]} W. Blaschke,
{\it Kreis und Kugel}. 

New York: Chelsea, 1949.

\item{[B\"oQu91]} 
J. B\"ohm and E. Quaisser
{\it Sph\"aroformen und symmetrische K\"orper,} 
Berlin: Akademie Verlag, 1991.

\item{[BoFe34]} 
T. Bonnesen and W. Fenchel
{\it Theorie der konvexen K\"orper,} 
Berlin: Springer, 1934.

\item{[BuZa88]} Yu. D. Burago and V. A. Zalgaller,
{\it Geometric Inequalities}, Grundlehren der mathematischen 
Wissenschaften {\bf 285}.Berlin: Springer, 1988.

\item{[CCG96]} 
S. Campi, A. Colesanti, and P. Gronchi,
Minimum problems for
volumes of constant bodies, in: 
P. Marcellini, G. Talenti, and E. Visintin, eds., 
{\it Partial Diferential Equations and Applications,} 

New York: M. Dekker, 1996, p. 43--55.

\item{[ChGr83]} G.D. Chakerian and H. Groemer,
Convex bodies of constant
width, in:  P.M. Gruber and J.\-M. Wills, eds.,
{\it Convexity and its applications,}
Basel: Birkh\"auser, 1983, 
p.  49--96.

\item{[Egg52]} H.G. Eggleston,
A proof of Blaschke's theorem on the Reuleaux triangle,
{\it Quart. J. Math. Oxford} ser. 2 {\bf 3} (1952)296--297.


\item{[Fey89]} R.P. Feynman,
{\it What do you care what
other people think?}.  New York: W.W.Norton, Bantam edition,
New York,1989, pp.166--168.

\item{[Fuj27-31]} M. Fujiwara,
Analytical proof of Blaschke's Theorem on the
curve of constant breadth with minimum area I and II,
{\it Proc.  Imp.  Acad. Japan} {\bf 3} (1927)307--309, and {\bf 7} (1931)300--30   2.

\item{[Gha96]} M. Ghandehari, An optimal control formulation 
of the Blaschke--Lebesgue Theorem,
{\it J.  Math.  Anal.  Appl.} {\bf 200} (1996)322--331.


\item{[Ka98]} B. Kawohl,
Was ist eine runde Sache? 

{\it GAMM Mitt. Ges. Angew. Math.
Mech.} {\bf 21} (1998)43--56.

\item{[Leb14]} H. Lebesgue, Sur le probl\`eme des isop\'erim\`etres et sur les d   omaines de largeur constante.
{\it Bull. Soc. Math. France}  (7) {\bf 4} (1921)67--96.

\item{[Leb21]} H. Lebesgue, Sur quelques questions des minimums, relatives
aux courbes orbiformes, et sur les rapports avec le calcul de variations.
{\it J.  Math.  Pure Appl.}  (7) {\bf 4} (1921)67--96.

%
\item{[Min67]} H. Minkowski, gesammelte Abhandlungen I and II.  New York: Chelse   a Publishing Co., 1967 (original, Leipzig, 1911).

\item{[M\"ul66]} K. M\"uller,
{\it Spherical Harmonics}, Springer Lecture Notes In Mathematics 
{\bf 17}.  Berlin: Springer, 1966.

\item{[Web94]} R.J. Webster,
{\it Convexity}.  Oxford: Oxford University Press, 1994.

\item{[Wei1884]} J. Weingarten, Festschrift der technischen Hochschule.  Berlin,  1  884.  (As cited in [Bla49].)

\end